\documentclass[12pt]{amsart}
\usepackage{amsmath,amsfonts,amssymb}

\def\ssum{\mathop{\sum\!\sum}}

\def\sumstar{\sideset{}{^*}\sum}

\def\le{\leqslant}

\def\ge{\geqslant}

\renewcommand{\mod}{\mathop{\rm{mod}}}

\def\cL{\mathcal{L}}
\def\cP{\mathcal{P}}
\def\cQ{\mathcal{Q}}
\def\cR{\mathcal{R}}

\numberwithin{equation}{section}

\theoremstyle{definition}

\begin{document}
\title{\bf Remarks on the Bombieri-Davenport \\ Large Sieve Inequalities}

\author[Friedlander]{John Friedlander$^*$}
\thanks{$^*$\ Supported in part by NSERC grant A5123}
\author[Iwaniec]{Henryk Iwaniec}

\maketitle

\dedicatory \quad\quad\quad\quad\quad\quad\quad\quad\quad\quad\quad\quad\quad\quad\quad{\sl to Enrico, just because!}
\medskip

{\bf Abstract:} Improvements of the Large Sieve for special sequences.

%The Large Sieve Inequality, although best possible in general, is subject to some improvements for some basic special sequences. We exemplify how ideas of Bombieri and Davenport can be exploited to this end.

\section{\bf Introduction: A bit of history}

It was just fifty years ago at the Vancouver ICM that Enrico Bombieri received
the Fields Medal. This sweet anniversary provokes us to refresh the memory of some of Enrico's early achievements in Analytic Number Theory. We choose to write about the Large Sieve Inequalities for Dirichlet characters, in part because they are fundamental to the theory of $L$-functions and in part because Enrico's role was so central to their development.

The modern version of the LSI asserts that
\begin{equation}\label{eq:1.1}
\sum_{q\le Q}\frac{q}{\varphi (q)}\sumstar_{\chi(\mod q)}\bigl|\sum_{M<n\le M+N}a_n\chi(n)\bigr|^2 \le (N+Q^2)  \sum_{M<n\le M+N}|a_n|^2 ,
\end{equation}
the asterisk restricting the summation to primitive characters $\chi (\mod q)$. This holds for all complex numbers $a_n$ supported on any interval of length $N\ge 1$. In numerous applications to prime numbers the
LSI serves as a substitute for the Riemann Hypothesis for $L(s,\chi)$. In fact, the LSI is found to be quite robust. Using~\eqref{eq:1.1}, one can recover territories beyond the threshold of the GRH.

The first ideas of the Large Sieve appear in the short paper of Ju. V. Linnik
[Li] followed by that of A. R\'enyi [Re]. Great innovations began in the 1960's with the
works of K. Roth [Ro] and E. Bombieri [B1]. Various conceptual interpretations
of the LSI have been gradually revealed, the common ground of all lying in the
super-orthogonality of primitive characters, the pseudo-characters of Selberg
[S1] included.

Many results, with $N+Q^2$ in~\eqref{eq:1.1} replaced by $c_1N+c_2Q^2$ or by
$\max (c_1N,c_2Q^2)$ for specific constants $c_1, c_2$, have been produced by
different authors, among them Davenport, Halberstam, Gallagher and Elliott. Ultimately,
the inequality~\eqref{eq:1.1}, best possible in general, was established independently by Montgomery and Vaughan [MV] and by Selberg [S2] by different arguments. 

For our considerations in this work, especially in Section 3, the value of $c_1=1$ is crucial while that of $c_2$ is less important. The first LSI of that strength in the literature is
the following.

{\bf Theorem 1.1:} (Bombieri-Davenport 1969) {\sl For all complex numbers $a_n$
supported on an interval of length $N\ge 1$ we have }
\begin{equation}\label{eq:1.2}
\sum_{q\le Q}\frac{q}{\varphi (q)}\sumstar_{\chi(\mod q)}\bigl|\sum_na_n\chi(n)\bigr|^2 \le ({\sqrt N} +Q)^2  \sum_n|a_n|^2 .
\end{equation}

\medskip

Because, in these notes, $Q$ will be much smaller than $\sqrt N$, the BD
inequality~\eqref{eq:1.2} is, for our purposes, as good as the MV-S inequality~\eqref{eq:1.1}. They are essentially sharp in general, however for special coefficients $a_n$ some improvements are possible. We have 

{\bf Theorem 1.2:} (Bombieri-Davenport 1969) {\sl Let $\cQ$ be a set of positive integers $\le Q$ which have no prime divisors in a set $\cP$. Then }
\begin{equation}\label{eq:1.3}
\sum_{q\in \cQ}\sum_{\chi(\mod q)}\frac{|\tau(\chi)|^2}{\varphi (q)}\bigl|\sum_n a_n\chi(n)\bigr|^2 \le ({\sqrt N} +Q)^2  \sum_n|a_n|^2 
\end{equation}
{\sl for all complex numbers $a_n$ supported on an interval of length $N\ge 1$ with $n$ being free of prime divisors in $\cP$.}

\medskip

Note that $\chi$ in~\eqref{eq:1.3} ranges over all characters in $\cQ$ not only the
primitive characters as in~\eqref{eq:1.2}. However, every $\chi (\mod q)$
with $1\le q\le Q$ is induced by a unique primitive character of conductor $q_1$ with
$q=q_1r\le Q$. Taking all these characters into account one obtains
\begin{equation}\label{eq:1.4}
\ssum_{\substack{qr\le Q\\(q,r)=1}}\frac{q\mu^2(r)}{\varphi (qr)}\sumstar_{\chi(\mod q)}|\sum_na_n\chi(n)|^2 \le ({\sqrt N} +Q)^2 \sum_n|a_n|^2 
\end{equation}
for all complex numbers $a_n$ supported on an interval of length $N\ge 1$ with
$n$ being free of prime divisors $\le Q$. Because
\begin{equation}\label{eq:1.5}
  \sum_{\substack{r\le X\\(r,q)=1}}\frac{\mu^2(r)}{\varphi (r)}
  \ge \frac{\varphi (q)}{q}\log X , \quad {\rm if}\quad  X\ge 1 ,
\end{equation}
we have 
\begin{equation}\label{eq:1.6}
\sum_{q\le Q}\bigl(\log\frac{Q}{q}\bigr)\sumstar_{\chi(\mod q)}\bigl|\sum_na_n\chi(n)\bigr|^2 \le ({\sqrt N} +Q)^2  \sum_n|a_n|^2 .
\end{equation}
See the arguments on page 21 of [BD] and a slightly stronger inequality in Theorem 8 of [B2], wherein the right hand side is the same as the right hand side of~\eqref{eq:1.1}.

Discarding all but the principal character in~\eqref{eq:1.6}
and taking $a_p=1$ for $M<p\le M+N$ and zero elsewhere, Bombieri and Davenport show that (choose
$Q=\sqrt N / \log N$)
\begin{equation*}
\pi (M+N) -\pi(M)\le \frac{2N}{\log N}\biggl(1+O\bigl(\frac{\log\log N}{\log N}\bigl)\biggl) ,
\end{equation*}
if $M>\sqrt N $. Recall that the linear sieve (whether Selberg or beta) produces the same interesting constant 2 and similar error terms. 

Our task in these notes is to exploit further the Bombieri-Davenport large sieve arguments, first in Section 2 for coefficients which have restricted support. Then, in Section 3, we show how an exceptional real character interacts with all of the others.
Our point of departure for these results is the following specific
inequality which occurs with a nice proof as Theorem 7A in [B2].

\medskip

{\bf Theorem 1.3:} {\sl For all complex numbers $a_n$ supported on
an interval of length $N\ge 1$ we have}
\begin{equation}\label{eq:1.7}
\ssum_{\substack{qr\le Q\\(q,r)=1}}\frac{q}{\varphi (qr)}\sumstar_{\chi(\mod q)}|\sum_na_n\chi(n)c_r(n)|^2 \le (N +Q^2) \sum_n|a_n|^2 ,
\end{equation}
{\sl where $c_r(n)$ is the Ramanujan sum}
\begin{equation}\label{eq:1.8}
c_r(n)=\sumstar_{u(\mod r)}e(un/r) = \sum_{d|(n,r)}d\mu(r/d) .
\end{equation}

\bigskip

\section{\bf The Large Sieve of Restricted Support}

Our first observation is that~\eqref{eq:1.7} implies the following.

\medskip

{\bf Proposition 2.1:} 
  {\sl Let $\cR$ be a set of numbers $r\le R$. Put}
  \begin{equation}\label{eq:2.1}
    \cL_q(\cR)= \frac{q}{\varphi (q)}\sum_{\substack{r\in \cR \\(r,q)=1}}\frac{\mu^2(r)}{\varphi (r)} ,
    \end{equation}
 \begin{equation}\label{eq:2.2}
\cL (Q,R) =\min_{q\le Q}\cL_q(\cR) . 
\end{equation} 
{\sl  Then}
\begin{equation}\label{eq:2.3}
 \sum_{q\le Q}\,\sumstar_{\chi(\mod q)}|\sum_{M<n\le M+N}a_n\chi(n)|^2 \le \frac{Q^2R^2+N}{\cL(Q,R)}  \sum_{M<n\le M+N}|a_n|^2 
 \end{equation} 
{\sl for sequences $a_n$ having $(n,r)=1$ for every $r\in \cR$.}

\medskip

{\bf Proof:}
  This follows from the fact that $c_r(n) =\mu(r)$ in this case.

\medskip

{\bf Example:}
Taking $\cR$ to be the set of all numbers $r\le R$ with $R\ge 2$, we get
\begin{equation}\label{eq:2.4}
\cL (Q,R) \ge \log R . 
\end{equation}

The number of primitive characters $\chi \mod q$ with $q\le Q$ is about $Q^2$
%asymptotically
%$$
%Q^2\prod_p\bigl(1-\frac{1}{p^2}-\frac{1}{p^3}\bigr) ,
%$$
so we have lost the factor $R^2$ in the first term of $Q^2+N$ in order to gain
$\log R$ in both terms. Although, at first glance this seems wasteful, it
could be only a minor complication in some important applications when the
character sums $\sum a_n\chi(n)$ are quite a bit longer than $Q^2$. 

\medskip

The previous example covers the situation where the $a_n$ are supported
on the primes in an interval. With more work and only half the success we 
consider what might be regarded the next natural example. 

\medskip

{\bf Proposition 2.2:} 
{\sl  Let $r(n)$ denote the number of representations of $n$ as the sum of two co-prime squares. We have}
\begin{equation}\label{eq:2.5}
 \sum_{q\le Q}\,\sumstar_{\chi(\mod q)}|\sum_{M<n\le M+N}r(n)a_n\chi(n)|^2 \le \frac{2N}{\sqrt{\log N/Q^2}} \sum_{M<n\le M+N}r^2(n)|a_n|^2 
 \end{equation} 
{\sl if $N\ge Q^2\exp(\alpha \log\log Q)^3$, where $\alpha$ is an
absolute constant.}

\medskip

To prove~\eqref{eq:2.5} using Proposition 2.1 we want to show that
$ \cL (Q,R) \gg (\log R)^{1/2}  $
for the relevant choice of $\cR$.
In preparation, we are first going to evaluate asymptotically
\begin{equation}\label{eq:2.6}
 S_q(x)=\sum_{\substack{n\le x\\(n,q)=1}}\nu(n)\tau(n)n^{-1} ,
\end{equation}
where $\nu(n)=1$ if $n$ is the product of distinct primes $\equiv 3 (\mod 4)$
and $\nu(n)=0$ otherwise. 

\medskip

{\bf Lemma 2.1:} {\sl For every $x\ge 1$ we have}
\begin{equation}\label{eq:2.7}
  S_q(x)=c\prod_{p|q_3}\bigl(1+\frac{2}{p}\bigr)^{-1} \log x
  +O\bigl( (1+\sum_{p|q_3}p^{-1}\log p)\prod_{p|q_3}(1+\frac{2}{p})\bigr) ,
  \end{equation}
    {\sl where the implied constant is absolute, where $q_3$ denotes the product of all distinct prime divisors of $q$ which are
      $\equiv 3 (\mod 4)$ and}
\begin{equation}\label{eq:2.8}
c=\frac{3}{\pi}\prod_{p\equiv 3 (\mod 4)}\bigl(1-\frac{2}{p(p+1)}\bigr) .
\end{equation}

\medskip

{\bf Proof:} We start from the generating Dirichlet series
\begin{equation}\label{eq:2.9}
  Z_q(s)=\sum_{(n,q)=1}\frac{\nu(n)\tau(n)}{n^s}=\prod_{p\nmid q_3}
  \bigl(1+ \frac{1-\chi_4(p)}{p^s}\bigr) ,
\end{equation}
where $\chi_4$ is the character of conductor 4. For $q=1$ this yields
\begin{equation}\label{eq:2.10}
\begin{aligned}
  Z(s) & = \frac{\zeta (s)}{L(s,\chi_4)}\prod_p\bigl(1-\frac{1}{p^s}\bigr)
  \bigl(1-\frac{\chi_4(p)}{p^s}\bigr)^{-1}\bigl(1+\frac{1-\chi_4(p)}{p^s}\bigr)\\
  & = \frac{\zeta (s)}{L(s,\chi_4)}\bigl(1-\frac{1}{4^s}\bigr)
  \prod_{p\equiv 3 (\mod 4)}\bigl( 1+\frac{1}{p^s}-\frac{2}{p^{2s}}\bigr)
 \bigl(1+\frac{1}{p^s}\bigr)^{-1}. 
  \end{aligned}
\end{equation}
Here the product over $p\equiv 3 (\mod 4)$ converges absolutely in ${\rm Re}(s)>\frac 12 $. Using the zero-free region for $L(s, \chi_4)$ (or elementary methods) we obtain
\begin{equation}\label{eq:2.11}
\sum_{n\le x}\nu(n)\tau(n)=cx +O\bigl(x(\log x)^{-2}\bigr).
\end{equation}
Note that $c$ is the residue of $Z(s)$ at $s=1$ and $L(1,\chi_4) =\pi/4$
(Leibniz formula). Hence, by partial summation  we get
\begin{equation}\label{eq:2.12}
S(x)= \sum_{n\le x}\nu(n)\tau(n)n^{-1}=c\log x +O(1) 
\end{equation}
which is the case of~\eqref{eq:2.7} for $q=1$.  In general we write
\begin{equation}\label{eq:2.13}
Z_q(s)=Z(s)F_q(s)
\end{equation}
where
\begin{equation*}
  F_q(s)= \prod_{p|q_3}\bigl( 1+\frac{2}{p^s}\bigr)^{-1}
  =\sum_{a=1}^{\infty}\frac{f(a)}{a^s}
\end{equation*}
and $f$ is multiplicative with $f(p^{\alpha})= (-2)^{\alpha}$ if $p|q_3$
and $f(p^{\alpha})= 0 $ if $p\nmid q_3$. Therefore, the coefficients of
$Z_q(s)$ are the convolution of the coefficients of $Z(s)$ with those of
$F_q(s)$. Hence our sum~\eqref{eq:2.6} is equal to (see~\eqref{eq:2.12})
\begin{equation*}
  \begin{aligned}
  S_q(x) &  =\ssum_{am\le x}f(a)\nu(m)\tau(m)(am)^{-1}
  =\sum_{a\le x}\frac{f(a)}{a}\bigl(c\log \frac{x}{a} + O(1)\bigr) \\
 &  = c\bigl(\sum_a\frac{f(a)}{a}\bigr)\log x
  +O\bigl(\sum_a\frac{|f(a)|}{a}\log 2a\bigr) .
   \end{aligned}
 \end{equation*}
The first sum above is equal to $F_q(1)$ which appears in~\eqref{eq:2.7}. The
second sum is bounded up to a constant by
\begin{equation*}
  \sum_a\frac{|f(a)|}{a}\bigl(1+\sum_{p|a}\log p\bigr)
  \ll \bigl(1+\sum_{p|q_3}\frac{\log p}{p}\bigr)
  \prod_{p|q_3}\bigl( 1 + \frac{2}{p}\bigr) .
\end{equation*}
This completes the proof of Lemma 2.1

\medskip

Estimating the error terms in Lemma 2.1 as follows:
\begin{equation*}
\sum_{p|q_3}p^{-1}\log p \ll \sum_{p<2\log q}p^{-1}\log p \ll \log\log 3q ,
\end{equation*}
  \begin{equation*}
    \prod_{p|q_3}\bigl(1+\frac{2}{p}\bigr)^2\ll \prod_{\substack{p<2\log q\\
        p\equiv 3(\mod 4)}}\bigl(1+\frac{1}{p}\bigr)^4 \ll (\log\log 3q)^2 ,
  \end{equation*}
  we may write~\eqref{eq:2.7} in the form
\begin{equation}\label{eq:2.14}
  S_q(x)=c\prod_{p|q_3}\bigl(1+\frac{2}{p}\bigr)^{-1}\bigl(\log x
  + O((\log\log 3q)^3)\bigr) .
\end{equation}
  Hence
\begin{equation}\label{eq:2.15}
  S_q(x)\gg \prod_{p|q_3}\bigl(1+\frac{2}{p}\bigr)^{-1}\log x
\end{equation}
if
\begin{equation}\label{eq:2.16}
x \ge \exp(\alpha\log\log 3q)^3 ,
\end{equation}
with some absolute constant $\alpha \ge 1$. Next, we remove
the divisor function $\tau(n)$ in~\eqref{eq:2.6} and estimate
\begin{equation}\label{eq:2.17}
T_q(x)= \sum_{\substack{n\le x \\ (n,q)=1}} \nu(n)n^{-1} .
\end{equation}
Clearly $T_q(x)^2\ge S_q(x)$, so~\eqref{eq:2.15} yields
\begin{equation}\label{eq:2.18}
T_q(x)\gg \prod_{p|q_3}\bigl(1+\frac{1}{p}\bigr)^{-1}(\log x)^{1/2}
\end{equation}
if $x$ satisfies~\eqref{eq:2.16}. This implies that~\eqref{eq:2.2}
satisfies
\begin{equation}\label{eq:2.19}
\cL (Q,R) \gg (\log R)^{1/2}
\end{equation}
with $R=\sqrt N Q^{-1}$, provided that $N>Q^2\exp(\alpha\log\log Q)^3$.
This completes the proof of Proposition 2.2. 
  
\bigskip

Next, we consider coefficients $a_n$ with $a_n\sqrt n$ small whenever $n$ has a relatively small prime divisor. One may think of $n$ as being almost a prime.

\medskip

{\bf Theorem 2.1}
{\sl Suppose that $a_n$ with $M < n < M+N$ satisfy }
\begin{equation}\label{eq:2.20}
  \sum_n|a_n|^2 \le 1 , \quad \sum_{n\equiv 0(\mod p)}|a_n|^2\le \frac{1}{p}
  \bigl(\frac{\log p}{\log N}\bigr)^2 .
\end{equation}
{\sl We have}
\begin{equation}\label{eq:2.21}
  \sum_{q\le Q}\, \sumstar_{\chi(\mod q)}|\sum_{M<n\le M+N}a_n\chi(n)|^2
  \le 24 N(\log N/Q^2)^{-1} 
\end{equation}
{\sl if } $8Q^2\le N$. 

\medskip

{\bf Proof:} Let $2\le R, \, Q^2R^3<N$. For the character sums $\sum a_n \chi(n)$ running over $n$ having no prime divisors $\le R$ we apply~\eqref{eq:2.3}
and~\eqref{eq:2.4}, obtaining the bound
\begin{equation*}
S_1 =\frac{2N}{\log R}\sum_n|a_n|^2\le \frac{2N}{\log R} .
\end{equation*}
  Then we estimate the character sums  $\sum a_n \chi(n)$ over $n$ having a prime divisor $\le R$ by 
\begin{equation*}
\sum_{p \le R}\bigl|\sum_{(n,P(p))=1}a_{np}\chi(n)\bigr| ,
\end{equation*}
where $P(z)$ denotes the product of all primes $<z$.
%Hence
This part of the sum~\eqref{eq:2.21} is bounded by
$S_2=\sum_{p_1}\sum_{p_2}S(p_1,p_2)$ with $p_1, p_2\le R$ and
\begin{equation*}
    \begin{aligned}
  S(p_1,p_2) = \sum_{q\le Q}\, \sumstar_{\chi (\mod q)} 
  |\sum_{(n,P(p_1))=1}a_{np_1}\chi(n)|  |\sum_{(n,P(p_2))=1}a_{np_2}\chi(n)|\\
   \le S(p_1,p_1)^{1/2}S(p_2,p_2)^{1/2} .
    \end{aligned}
\end{equation*}

Applying~\eqref{eq:2.3} we get
\begin{equation*}
  S(p,p) \le \frac{Q^2p^2+N/p}{\log p}\sum_n|a_{np}|^2
  \le \frac{2N\log p}{(p\log N)^2} .
\end{equation*}
Hence,
\begin{equation*}
  S_2 \le \frac{2N}{(\log N)^2}\bigl(\sum_{p\le R}\frac{1}{p}
  \sqrt{\log p}\bigr)^2\le 18 N(\log N)^{-2}\log R\le 2N(\log R)^{-1} .
\end{equation*}
This completes the proof of Theorem 2.1 on choosing $R=(N/Q^2)^{1/3}$.

\medskip

{\bf Remark:}
We could assume the somewhat stronger bound
\begin{equation}\label{eq:2.22}
a_n\ll \frac{\rho(n)}{\sqrt n}, \quad 1\le n\le N
\end{equation}
with
\begin{equation}\label{eq:2.23}
\rho(n)\ = \prod_{p|n}\rho (p),  \quad \rho(p)=\frac{\log p}{\log N} .
\end{equation}
These coefficients satisfy~\eqref{eq:2.20} up to a constant. Indeed, we have

\begin{equation*}
\sum_{n\le N}\frac{\rho(n)^2}{n}\ll \prod_{p\le N}\bigl(1+\frac{\rho(p)^2}{p}\bigr)\le \exp\biggl(\sum_{p\le N}\frac{\rho(p)^2}{p}\biggr) ,
\end{equation*}
\begin{equation*}
\sum_{p\le N}\frac{\rho(p)^2}{p}\ll \sum_{p\le N}\frac{1}{p}\biggl(\frac{\log p}{\log N}\biggr)^2\ll 1 .
\end{equation*}
The above estimates verify both conditions in~\eqref{eq:2.20}.

\medskip 

{\bf Example:} The coefficients
\begin{equation*}
a_n=n^{-1/2}\Lambda_k(n)(\log N)^{-k} , \quad {\rm with} \quad 1\le n \le N , 
\end{equation*}
satisfy~\eqref{eq:2.20}. Here $\Lambda_k $ is the generalized von-Mangoldt
function $\Lambda_k= \mu * \log^k $. Apply induction in $k$ using
$\Lambda_{k+1}=\Lambda_k\cdot \log +\Lambda_k * \Lambda$. 

\bigskip

\section{\bf The Large Sieve with an Exceptional Character} 

Recall that taking only the one term of~\eqref{eq:1.6} coming from the principal character over primes in a segment Bombieri and Davenport obtained a strong Brun-Titchmarsh upper bound for the primes in an interval. 
One may additionally extract the contribution of a second real character, say  $\chi_D$ of conductor $D$, which is 
\begin{equation*}
\bigl(\log \frac{Q}{D}\bigr) \biggl(\sum_{M<p\le M+N}\chi_D(p)\biggr)^2 . 
\end{equation*}
In this section we study some of the effects of doing so. 

We put $\lambda(p)= 1+\chi_D(p)$ and write
\begin{equation*}
  \bigl(\sum_p\chi_D(p)\bigr)^2 = \bigl(-\sum_p1 + \sum_p\lambda(p)\bigr)^2
  \ge  \bigl(\sum_p1\bigr)^2 -2\bigl(\sum_p1\bigr)\bigl(\sum_p\lambda(p)\bigr) .
\end{equation*}
Hence, Theorem 8 of [B2] yields 

\begin{equation}\label{eq:3.1}
  \begin{aligned}
    \sum_{1<q\le Q}(\log\frac{Q}{q}) &   \sumstar_{\substack{\chi \mod q\\\chi\neq \chi_D}}  |\sum_p\chi(p)|^2  \le(Q^2+N)(\sum_p 1) \\ & -(\log\frac{Q^2}{D})(\sum_p1)^2 +2(\log\frac{Q}{D})(\sum_p1)(\sum_p\lambda(p)) .
  \end{aligned}
\end{equation}

Note that the last part of~\eqref{eq:3.1} is negligible if $\chi_D(p)= -1$ very frequently. We are going to use the above arrangement for the coefficients
in~\eqref{eq:1.6} given by 
\begin{equation}\label{eq:3.2}
a_n=\Lambda(n)f(n/N)
  \end{equation}
where $f(u)$ is a nice function supported on $0\le u\le 1$ and with $0\le f(u)\le 1$. Let $3\le Q\le\sqrt N$. We get
\begin{equation}\label{eq:3.3}
  \begin{aligned}
    \sum_{q\le Q}(\log\frac{Q}{q})\sumstar_{\chi(\mod q)}& |\sum_n\chi(n)\Lambda(n)f(n/N)|^2 \\ & \le (Q^2+N)   \sum_n|a_n|^2 +O\bigl(N^{7/4}\log N\bigr) ,
 \end{aligned}
\end{equation}
where the error term on the right side comes from an estimation of the contribution of prime powers $ p^2, p^3, \ldots \le N$ and primes $p\le Q$.

For $q=1$ we have the principal character whose contribution is 
\begin{equation}\label{eq:3.4}
(\log Q)(\sum_n a_n)^2 . 
\end{equation}

Let $\chi_D$ be a primitive real character of conductor $D\le Q$. Its contribution is  
\begin{equation}\label{eq:3.5}
\bigl(\log \frac{Q}{D}\bigr)\bigl(\sum_n \chi_D(n)a_n\bigr)^2. 
\end{equation}

Put $\lambda=1*\chi_D$ so $\chi_D=\mu * \lambda$. Note that
$\lambda(p)=1+\chi_D(p)$. We think of $\lambda$ as being lacunary. We write
\begin{equation}\label{eq:3.6}
\sum_n \chi_D(n)a_n = -\sum_na_n +\sum_n\rho(n)a_n 
\end{equation}
where $\rho(n)=1+\chi_D(n)$. Hence
\begin{equation}\label{eq:3.7}
(\sum_n \chi_D(n)a_n)^2 \ge (\sum_na_n)^2- 2(\sum_na_n)(\sum_n\rho(n)a_n). 
\end{equation}
Summing up the contributions of $\chi_0, \chi_D$, we obtain from~\eqref{eq:3.3} 

\begin{equation}\label{eq:3.8}
  \begin{aligned}
    \sum_{1<q\le Q} (\log\frac{Q}{q}) & \sumstar_{\substack{\chi(\mod q)\\ \chi\neq \chi_D}}|\sum_n\chi(n)\Lambda(n)f(n/N)|^2
  \\ &  \le (Q^2+N) \sum_na_n^2
     -\bigl(\log\frac{Q^2}{D}\bigr)(\sum_na_n)^2
  \\ &  +2\bigl(\log\frac{Q}{D}\bigr)
    (\sum_na_n)(\sum_n\lambda(n)a_n) 
    +O\bigl(N^{7/4}\log N\bigr) .
  \end{aligned}
  \end{equation}

Note that we have here replaced $\rho(n)$ by $\lambda(n)$ by estimating trivially the terms at higher prime powers. Next, we estimate as follows.
\begin{equation*}
  \begin{aligned}
  \sum_n \lambda(n)a_n & \le \sum_{n\le N}\lambda(n)\log N \\ & =
  L(1,\chi_D)N\log N +O(N^{1/2}D^{1/4}\log N) ,
   \end{aligned}
\end{equation*}
where $L(1,\chi_D)$ is the residue of $\zeta(s)L(s,\chi_D)$
at $s=1$.

Moreover, we have $ \sum_n a_n \ll N$, but we need
more precise estimations for the leading terms.
Using the prime number theorem we get
\begin{equation}\label{eq:3.9}
  \sum_n a_n= A_1N +O(N/\log N) , \quad  A_1= \int_0^1f(u)du ;
\end{equation}
\begin{equation}\label{eq:3.10}
  \sum_n a_n^2 = A_2N\log N +O(N), \quad A_2= \int_0^1f^2(u)du .
\end{equation}
%\begin{equation}\label{eq:3.11}
%  A_1= \int_0^1f(u)du , \quad\quad  A_2= \int_0^1f^2(u)du .
%\end{equation}
Introducing these estimates into~\eqref{eq:3.8} and choosing $Q={\sqrt N}/\log N$, we obtain the following.

\medskip

{\bf Lemma 3.1:}
{\sl   Let $3\le D \le Q={\sqrt N}/\log N$. We have}
\begin{equation}\label{eq:3.12}
  \begin{aligned}
    \sum_{1<q\le Q}(\log\frac{Q}{q})
   &  \sumstar_{\substack{\chi(\mod q)\\ \chi\neq \chi_D}}
    |\sum_n\chi(n)\Lambda(n)f(n/N)|^2 \\ &
    \le (A_2\log N -A_1^2\log N/D)N^2 \\ &
    +2A_1(\log Q/D)N^2L(1,\chi_D)\log N
    +O(N^2) .
     \end{aligned}
     \end{equation}

\medskip

Choosing $f(u)=1$ we find $A_1=A_2=1$. Hence, 

\medskip

{\bf Proposition 3.1:}
{\sl  Let $3\le D \le Q={\sqrt N}/\log N$. We have}
\begin{equation}\label{eq:3.13}
  \begin{aligned}
    \sum_{1<q\le Q}(\log\frac{Q}{q})
   &  \sumstar_{\substack{\chi(\mod q)\\ \chi\neq \chi_D}}
    |\sum_n\chi(n)\Lambda(n)|^2 \\ &
    \le N^2\bigl(\log cD +L(1,\chi_D)(\log N)^2\bigr) ,
     \end{aligned}
     \end{equation}
{\sl where $c\ge 1$ is an absolute constant.}

\medskip

{\bf Proposition 3.2:} {\sl Let $D$ be larger than a suitable absolute constant. Suppose}
\begin{equation}\label{eq:3.14}
L(1,\chi_D)\log D \le {\varepsilon}^5.
  \end{equation}
{\sl Then, for every primitive character $\chi (\mod q), \chi \ne \chi_0, \chi_D$ and every $N\ge q^4$ with $D^{1/\varepsilon^2}<N<D^{1/\varepsilon^3}$ we have}
\begin{equation}\label{eq:3.15}
|\sum_{n\le N}\chi (n) \Lambda (n)| \le 3 \varepsilon N .
\end{equation}

\medskip

{\bf Proof:} On the left side of~\eqref{eq:3.13} we have
$Q/q \ge N^{1/4}/\log N$. Hence,
\begin{equation*}
  |\sum_{n\le N}\chi (n) \Lambda (n)|^2 < \frac{9}{2}
  \bigl(\frac{\log D}{\log N} +\varepsilon^5\frac{\log N}{\log D}\bigr)N^2
< (3\varepsilon N)^2 
\end{equation*}
since $c$ is an absolute constant, $D, N$ are sufficiently large and $4<9/2$.

%%%%%%%%%%%%%%%%%%%%%%%%%%%%%%

\medskip 
Department of Mathematics, University of Toronto

Toronto, Ontario M5S 2E4, Canada  \quad (frdlndr@math.toronto.edu)

\medskip

Department of Mathematics, Rutgers University

Piscataway, NJ 08903, USA  \quad  (iwaniec@comcast.net)

\end{document}